\documentclass[psamsfont]{amsart}
\usepackage[utf8]{inputenc}
\usepackage{graphicx}
\usepackage[dvips]{epsfig}
\usepackage{pinlabel}
\usepackage{amsmath}
\usepackage{amsfonts}
\usepackage{latexsym}
\usepackage{amssymb}
\usepackage[usenames,dvipsnames]{color}
\usepackage{amsthm}
\usepackage[all]{xypic}
\usepackage{enumitem}
\usepackage[breaklinks=true]{hyperref}

\input xy

\xyoption{all}

\author{Emmy Murphy}
\address{Northwestern University, United States.}
\email{e\_murphy@math.northwestern.edu}

\theoremstyle{plain}
\newtheorem{thm}{Theorem}[section]
\newtheorem{lemma}[thm]{Lemma}
\newtheorem{prop}[thm]{Proposition}
\newtheorem{coro}[thm]{Corollary}

\theoremstyle{definition}

\newtheorem{definition}[thm]{Definition}
\newtheorem{remark}[thm]{Remark}
\newtheorem{remarks}[thm]{Remarks}

\theoremstyle{remark}

\numberwithin{equation}{section}

\newcommand{\R}{\mathbb{R}}
\newcommand{\Z}{\mathbb{Z}}
\newcommand{\C}{\mathbb{C}}

\newcommand{\J}{\mathcal{J}}
\newcommand{\p}{\varphi}
\newcommand{\e}{\varepsilon}
\newcommand{\dd}{\partial}

\newcommand{\es}{\varnothing}
\newcommand{\op}{\operatorname}
\newcommand{\sse}{\subseteq}

\newcommand{\x}{\times}
\newcommand{\sm}{\setminus}

\newcommand{\wt}{\widetilde}
\newcommand{\wh}{\widehat}
\newcommand{\ol}{\overline}

\newcommand{\std}{\text{std}}

\begin{document}

\title{Loose Legendrian embeddings in high dimensional contact manifolds}

\begin{abstract}
We give an $h$--principle type result for a class of Legendrian embeddings in contact manifolds of dimension at least $5$. These Legendrians, referred to as loose, have trivial pseudo-holomorphic invariants. We demonstrate they are classified up to Legendrian isotopy by their smooth isotopy class equipped with an almost complex framing. This result is inherently high dimensional: analogous results in dimension $3$ are false.
\end{abstract}

\maketitle

\section{Introduction}\label{sec:intro}

Let $(Y, \xi)$ be a contact manifold, that is a smooth manifold $Y$ with a hyperplane field $\xi$ which is maximally non-integrable. Throughout the paper we let $2n+1$ be the dimension of $Y$, which we will always assume is at least $5$. A \emph{Legendrian embedding} is a map $f: \Lambda \to Y$ so that $df(T\Lambda) \sse \xi$, where $\Lambda$ is a smooth $n$--manifold. We are interested in the general question of when two Legendrian embeddings $f_0, f_1: \Lambda \to Y$ are isotopic through Legendrian embeddings.

Since the codimension of a Legendrian embedding is $n+1 > 3$, asking whether Legendrians are \emph{smoothly} isotopic is mostly solved with the $h$--cobordism theorem, for example, when $Y$ is simply connected. However, the Legendrian isotopy question is known to be very intricate \cite{EES, BST}. Indeed, for $n\geq 2$, there is not even a conjectural classification of any kind, even for any fixed $(Y, \xi)$ and fixed topology of $\Lambda$. The simplest possible invariant one can develop beyond smooth isotopy type is the \emph{formal Legendrian} isotopy type.

\begin{definition}
Let $\Lambda$ be a smooth $n$--manifold, and $(Y, \xi)$ a contact $(2n+1)$--manifold. A \emph{formal Legendrian embedding} is a pair $(f, F_s)$, where $f: \Lambda \to Y$ is a smooth embedding, and $F_s: T\Lambda \to TY$ is a homotopy of bundle maps covering $f$, so that:
\begin{itemize}
\item $F_0 = df$,
\item $F_s$ is fiberwise injective for all $s \in [0,1]$, and
\item the image of $F_1$ is contained in $\xi$, and furthermore its image is Lagrangian with respect to the linear conformal symplectic structure on $\xi$.
\end{itemize}
We note that every Legendrian embedding is a formal Legendrian embedding, by letting $F_s = df$ for all $s$. Furthermore, any formal Legendrian embedding $(f, F_s)$ where $F_s$ is constant in $s$ is necessarily a genuine Legendrian embedding.
\end{definition}

Since the space of Legendrian embeddings is contained in the space of formal Legendrian embeddings, we can say that two Legendrian embeddings $f_0, f_1: \Lambda \to (Y, \xi)$ are \emph{formally isotopic} if there is a path of formal Legendrian embeddings interpolating between them. The advantage of working with formal Legendrians is that they are purely algebro-geometric objects: the forgetful map taking a formal Legendrian embedding $(f, F_s)$ to its smooth embedding $f$ is a Serre fibration (onto its image, it is not a surjection). Thus asking whether two Legendrian embeddings are formally Legendrian isotopic is first asking if they are smoothly isotopic, and then seeing whether that smooth isotopy extends to a formal Legendrian isotopy, which ultimately boils down to computing some classes using obstruction theory in the homotopy groups of frame bundles. The same is true for asking higher parametric questions about the space of formal Legendrian embeddings. For a familiar example, two Legendrian knots $f_0, f_1: S^1 \to \R^3_\std$ are formally Legendrian isotopic if and only if they have the same rotation and Thurston-Bennequin numbers. We do a similar computation for higher dimensional Legendrians in the Appendix \ref{appx}.

The second important concept is that of a \emph{loose} Legendrian. We delay the complete definition until Section \ref{sec:loose}, but here we note the following. There is a universal fixed model Legendrian $\Lambda_\ell \sse B^{2n+1}_\std$, called a \emph{loose chart}, where $\Lambda$ is diffeomorphic to a properly embedded disk, and $B^{2n+1}_\std \sse (\R^{2n+1}, \xi_\std)$ is the open unit ball and $\xi_\std = \ker(dz - \sum_{i=1}^n y_idx_i)$. Then a connected Legendrian submanifold $\Lambda \sse (Y, \xi)$ is \emph{loose} if there is an open set $U \sse Y$ so that the pair $(U, U \cap \Lambda)$ is contactomorphic to the pair $(B^{2n+1}_\std, \Lambda_\ell)$. We also say a Legendrian embedding is loose if its image is loose.

\begin{thm}\label{thm:main iso}
Suppose $n \geq 2$, and let $f_0, f_1: \Lambda \to (Y, \xi)$ be two loose Legendrian embeddings, which are formally isotopic. Then they are Legendrian isotopic.
\end{thm}

We also note that every formal Legendrian isotopy class contains loose Legendrian embeddings. In fact, in any fixed formal Legendrian isotopy class, the space of loose Legendrians is $C^0$ dense. This is almost immediate from the definition: in any small neighborhood of a point on a Legendrian we can cut out a small disk and replace it with a Legendrian disk containing a loose chart. 

Theorem \ref{thm:main iso} completely classifies loose Legendrian embeddings up to isotopy. In fact we have a more general theorem which gives us an understanding of higher parametric families of loose Legendrians. Fix a smooth manifold $\Lambda$ and an open disk $D^n \sse \Lambda$. For a contact manifold $(Y, \xi)$, fix an open set $U \sse Y$, a contactomorphism $U \cong B^{2n+1}_\std$, and a parametrization $\p: D^n \to B^{2n+1}_\std$ of $\Lambda_\ell$.  Let $\mathcal{L}^{\text{form}}_\ell(\Lambda, U)$ be the space of all formal Legendrian embeddings $(f, F_s)$, so that $f^{-1}(U) = D^n$, $(f, F_s)$ is a genuine Legendrian embedding on $D^n$, and $f|_{D^n} = \p$ with respect to the chosen contactomorphism $U \cong B^{2n+1}_\std$. Briefly, $\mathcal{L}^{\text{form}}_\ell(\Lambda, U)$ is the space of all formal Legendrian embeddings with a fixed loose chart. \footnote{Technically, the data of $D^n$, $\p$, and the contactomorphism between $U$ and $B^{2n+1}_\std$ should all be included in the notation defining $\mathcal{L}^{\text{form}}_\ell(\Lambda, U)$, as all of this data is fixed.}

\begin{thm}\label{thm:main}
Fix $k>0$ and $n \geq 2$, and for $t \in D^k$ let $(f_t, F_{s,t})$ be a smooth family in $\mathcal{L}^{\text{form}}_\ell(\Lambda, U)$, so that $(f_t, F_{s,t})$ is a genuine Legendrian embedding for all $t \in \dd D^k$. Then the family $(f_t, F_{s,t})$ is isotopic through formal Legendrian embeddings, rel $\dd D^k$ (though not necessarily in $\mathcal{L}^{\text{form}}_\ell(\Lambda, U)$), to a family of genuine Legendrian embeddings.
\end{thm}

We now explain the strategy of the proof, and the structure of the paper. There are essentially three main ingredients. First, we show that any family of formal Legendrian embeddings can be isotoped to a family which is ``graphical'': we define this concept and prove it in Section \ref{sec:graphical}. This is a fairly weak result, it does not really make the formal embeddings ``more Legendrian'', rather it just gives us global $1$--jet coordinates to work in. Using this result, we can reinterpret the question of isotoping a formal Legendrian to a genuine Legendrian as a question of approximating an arbitrary $1$--jet on $\Lambda$ by a holonomic $1$--jet of a function. This cannot be done, but it can be done with mild hypersurface singularities, using the concept of wrinkled embeddings from \cite{EM}. This is the content of Section \ref{sec:wrinkled}. Reinterpreting things again in the Legendrian world, we see that any family of formal Legendrian embeddings can be isotoped to a family of ``wrinkled Legendrians''. The question now remains of how to resolve these singularities, particularly in a way that is smooth in parametric families. This cannot be done in general, but it gives rise to the definition of a loose chart: a smooth model which can approximate the wrinkled singularities. This is explained in Section \ref{sec:loose}, where loose charts are defined and where the proof of Theorem \ref{thm:main} is completed. Finally in Section \ref{sec:conclusion} we explain some basic results, including the proof of Theorem \ref{thm:main iso} and the relationship to Legendrian contact homology and Lagrangian fillability. For reference, in Appendix \ref{appx} we compute the formal Legendrian isotopy classes in a fixed smooth isotopy class.

Since the original preprint of this paper appeared, a number of further results related to loose Legendrians have been discovered. We do not intend to survey these results here, but we give a brief bibliography. A significant source of applications has been to the geometry of Weinstein manifolds, see \cite{flex, CE} and Lagrangians inside them \cite{EGL}. Subsequent analysis of their boundary dynamics have given related results about fillings of contact manifolds \cite{oleg}. Loose Legendrians are also known to be closely tied to the theory of high dimensional overtwisted contact manifolds, see \cite{BEM, CMP, H}. They are also useful in the construction of Lagrangian embeddings \cite{caps, YETI, toru, GDR}. Their orderability properties were studied in \cite{Liu, PPP}, and relations with norms on contactomorphism groups were explored in \cite{looseOB}.

\subsection*{Acknowledgements}
The author thanks Y.~Eliashberg for many discussion during the original creation of this paper. The author is additionally grateful to many people for discussions which helped clarify the exposition and concepts, including R.~Casals, K.~Cieliebak, V.~Colin, P.~Massot, K.~Siegel, and M.~Sullivan.

\section{Graphical submanifolds}\label{sec:graphical}

Recall that, for any smooth manifold $\Lambda$, the first jet space of $\Lambda$ is the contact manifold $\J^1\Lambda = (\R \x T^*\Lambda, \ker(dz - \lambda_\std))$, where $\lambda_\std \in \Omega^1T^*\Lambda$ is the tautological $1$--form. This space is equipped with the natural projection $\pi: \J^1\Lambda \to \\R \x \Lambda$, which is called the \emph{front projection}. In this section we will assume $(Y, \xi)$ is coorientable.\footnote{All results are true in the non-coorientable case by replacing $\J^1\Lambda$ with the $1$--jet bundle of another real line bundle over $\Lambda$; the proofs are identical.} The main result of this section states that the neighborhood of any formal Legendrian embedding can be globally modeled by an open set in $\J^1\Lambda$, after a possible formal Legendrian isotopy. The author is particularly grateful to K.~Cieliebak and Y.~Eliashberg for sharing an early manuscript of an upcoming book, where the same proposition was proved.

\begin{prop}\label{prop:graphical}
Let $f: \Lambda \to (Y, \xi)$ be a formal Legendrian embedding, covered by maps $F_s: T\Lambda \to TY$. Then, after a smooth isotopy from $f$ to $\wt f$, we can choose an open set $U \sse Y$ containing $\wt f(\Lambda)$ and a map $\p: U \to \J^1\Lambda$ which is a contactmorphism onto its image, so that $\pi \circ \p \circ \wt f$ is the identity map on $\Lambda$.

This can also be done smoothly in families: if $(f_t, F_{s,t})$ is a family of formal Legendrian embeddings for $t \in D^k$, then we obtain an isotopic family of maps $\wt f_t$, an open set $U \sse Y \x D^k$, and a smooth family of maps $\p_t: U \cap Y \x \{t\} \to \J^1\Lambda$, so that each is a contactomorphism onto its image and $\pi \circ \p_t \circ \wt f_t$ is the identity for all $t$. Additionally, if $f_t$ is Legendrian on some closed set $A \sse \Lambda \x D^k$, then we can take $\wt f_t = f_t$ on this set, so that $\p_t$ maps $\wt f_t(A \cap \{t\})$ to the zero section in $\J^1\Lambda$.
\end{prop}

We start by proving that, after isotopy, we can make a family of formal Legendrian embeddings into a family of $\e$--Legendrian embeddings. Recall that a smooth embedding $f: \Lambda \to (Y, \xi)$ is called an $\e$--Legendrian embedding if at every point $x \in \Lambda$ there is a Legendrian plane $P \sse \xi_{f(x)}$ so that the angle between $P$ and $df(T\Lambda)$ is less than $\e$, with respect to some given metric on $Y$. We fix an $\e>0$, sufficiently small so that in each fiber, the space of $\e$--Legendrian planes deformation retracts onto the space of Legendrian planes. Notice that in a fixed fiber $TY$, the space of $\e$--Legendrian planes is just the $\e$--neighborhood of the Legendrian Grassmannian, as a subset of the $n$--plane Grassmannian in $TY$.

We denote $\op{Gr}_k(Y)$ to be the space of $k$--planes in $TY$, i.e. the fiber bundle over $Y$ whose fiber over $y \in Y$ is $\op{Gr}_k(TY_y)$.

\begin{prop}\label{prop:eLeg}
Let $(f_t, F_{s,t})$ be a family of formal Legendrian embeddings for $t \in D^k$. Then we can isotope the family $\{f_t\}$ to a family of $\e$--Legendrian embeddings $\wt f_t$. If $(f_t, F_{s,t})$ is already Legendrian on some closed set $A \sse \Lambda \x D^k$, then we can take $\wt f_t = f_t$ on $A$.
\end{prop}

\begin{proof}
This follows from the convex integration theorem for directed embeddings \cite[Theorem 19.4,2]{EM}, where the set $A \sse \op{Gr}_nY$ in the theorem is the set of $\e$--Legendrian planes. Clearly this set is open, it remains to show that it is affine ample, as defined in the \cite[Section 19.4]{EM}. 

Fix a point $y \in Y$. First, assume that $S \sse TY_y$ is an isotropic $(n-1)$--plane in $\xi_y$. Then the set of all vectors $v$ in $TY_y$ so that $S \oplus \R\cdot v$ is a Legendrian plane is simply $S^{d\alpha \perp} \sse \xi_y$. Here $d\alpha \perp$ is symplectic orthogonal complementation with respect to the linear, conformally defined symplectic structure on $\xi$. Let $S^\perp$, be the Riemannian orthogonal complement of $S$ in $TY$, and let $\pi: S^\perp \to S^{d\alpha \perp} \cap S^\perp$ be the orthogonal projection. Notice that $S^{d\alpha \perp} \cap S^\perp$ is a $2$-plane, since $S \sse S^{d\alpha \perp}$. Then for $v \in S^\perp$, as long the angle between $v$ and $\pi(v)$ is less than $\e$, the angel between $S \oplus \R\cdot v$ and $S \oplus \R \cdot \pi(v)$ is less than $\e$ as well, so the plane $S \oplus \R\cdot v$ is $\e$--Legendrian. Therefore the set of vectors $v \in S^\perp$ making $S \oplus \R\cdot v$ an $\e$--Legendrian contains the solid $\e$--cone around $S^{d\alpha \perp} \cap S^\perp$. (The $\e$--cone around a linear subspace is the set of vectors whose angle with the subspace is less than $\e$.) Let $C \sse S^\perp$ be this $\e$--cone around $S^{d\alpha \perp} \cap S^\perp$. Notice that the convex hull of $C$ is the entire space. And furthermore, for any affine hyperplane $H' \sse S^\perp$, the convex hull of $C \cap H'$ is all of $H'$, since $S^{d\alpha \perp} \cap S^\perp$ is $2$--dimensional.

Working now in $TY$ rather than $S^\perp$, we see that the set of vectors $v \in TY$ making $S \oplus \R \cdot v$ an $\e$--Legendrian contains the set $C \x S \sse S^perp \x S = TY$, since $S \oplus \R \cdot v_1 = S \oplus \R \cdot v_2$ whenever $v_1 - v_2 \in S$. And for any affine hyperplane $H' \sse TY$ parallel to $S$, the convex hull of $C \x S \cap H'$ is equal to $H'$. This verifies affine ampleness, in the case where $S$ is isotropic. If $S$ is not isotropic, we can ask whether it makes an angle less than $\e$ with some isotropic $(n-1)$--plane $\wt S$. If so, we can repeat the same proof above replacing $S$ with $\wt S$: the vectors $v \in S^\perp$ so that $S \oplus \R\cdot v$ makes angle less than $\e$ with $\wt S \oplus \R\cdot \pi(v)$ contains an $\e$--cone around $\wt S^{d\alpha \perp} \cap S^\perp$. Thus we remain affine ample in this case as well. Finally, if $S$ does not make an angle less than $\e$ with any isotropic $(n-1)$--plane, then the set of vectors $v$ making $S \oplus \R \cot v$ an $\e$--Legendrian is empty, which is ample as well. This completes the proof that the $\e$--Legendrian condition is affine ample, allowing us to use \cite[Theorem 19.4,2]{EM}.
\end{proof}

Constructing neighborhoods defined in Proposition \ref{prop:graphical} essentially comes down to finding a suitable Legendrian foliation representing the fiber coordinates in $\J^1\Lambda$. We define a \emph{Legendrian submersion} of a contact manifold $(Y, \xi)$ to be a submersion $s: Y \to N$ for some smooth $(n+1)$--manifold $N$, so that $\ker ds \sse \xi$ at every point in $Y$. Note that on $\J^1\Lambda$, the projection map $\J^1\Lambda \to \R \x \Lambda$, which is called the \emph{front projection}, is a Legendrian submersion. This is the model we are attempting to build, thus an $h$--principle for Legendrian submersions will be a useful ingredient.

\begin{prop}\label{prop:micro}
The differential relation defining Legendrian submersions in microflexible, in the sense of \cite[Chapter 13]{EM}. 
\end{prop}

\begin{proof}
This is essentially a form of contact stability. Smoothly, any compactly supported deformation of submersions is induced by an isotopy of the domain. Thus a deformation of Legendrian submersions is induced by a smooth isotopy. If the original submersion is locally given by $(x, y, z) \mapsto (x, z)$ with $\xi = \ker(dz - \sum y_i dx_i)$ (and locally every Legendrian submersion is this), then the new contact structure obtained by pulling back the standard contact form by the smooth isotopy is of the form $a_0(x,y,z)dz - \sum a_i(x,y,z)dx_i$ for some functions $a_i$, since the kernel is Legendrian. We can furthermore assume that $a_0$ is non-vanishing if the deformation is $C^\infty$ small, and then the contact condition implies that $\det\left[\frac{\dd}{\dd y_j}\frac{a_i}{a_0}\right]$ is non-vanishing. Thus, we can perform a further smooth isotopy, tangent to the kernel of the submersion, so that the result is a contact isotopy. Thus, any small deformation of Legendrian submersions is induced by a contact isotopy. 

To prove microflexibility then, we just note that small deformations of Legendrian submersions are induced by contact isotopies, and contact isotopies on open sets can be extended by cutting off the Hamiltonian. In particular this applies to the $\theta_k$ pairs defining microflexibility.
\end{proof}

\begin{lemma}\label{lem:transverse}
Let $f_t: \Lambda \to Y$ be a family of $\e$--Legendrian embeddings for $t \in D^k$. Then we can find a family of injective bundle maps $P_t: T^*\Lambda \oplus \R \to f_t^*TY$, so that at each point: 
\begin{itemize}
\item $P_t(T^*\Lambda \oplus \R)$ is transverse to $T\Lambda \sse f_t^*TY$, 
\item $P_t(T^*\Lambda \oplus \{0\})$ is a Legendrian plane in $f_t^*\xi$, and
\item $P_t(\{0\} \oplus \R)$ is transverse to $f_t^*\xi$.
\end{itemize}
\end{lemma}

\begin{proof}
For definitiveness, choose a contact form and a compatible almost complex structure. For each $x \in \Lambda$, let $L^x_t \sse \xi_{f_t(x)}$ be a Legendrian plane which makes an angle less than $\e$ with $df_t(T\Lambda_x)$. Since there is a deformation retraction of $\e$--Legendrian planes onto Legendrian planes this choose can be made smoothly with respect to $x$ and $t$. Let $\pi^x_t: df_t(T\Lambda_x) \to L^x_t$ be the orthogonal projection, which is an isomorphism. We then define $P_t|_{T^*\Lambda \oplus \{0\}} = J \circ \pi^x_t \circ df_t$, since $J(L^t_x)$ is orthogonal to $L^t_x$ this is transverse to $df_t(T\Lambda_x)$. Define $P_t|_{\{0\} \oplus \R}$ to be scalar multiplication by the Reeb vector field, with is orthogonal to $\xi$ and therefore transverse to $df_t(T\Lambda_x)$.
\end{proof}

\begin{proof}[Proof of Proposition \ref{prop:graphical}:]
We begin by replacing our family of formal Legendrian embeddings by a family of $\e$--Legendrian embeddings, we continue to use the notation $f_t: \Lambda \to Y$ for this latter family. Let $N = \Lambda \x (0,1)$, and let $U_t \sse Y$ be a continuous family of tubular neighborhoods containing $f_t(\Lambda)$. First, we construct a formal Legendrian submersion from $U_t$ to $N$: by definition this is a continuous map $U_t \to N$ and a surjective bundle map $G_t: TU_t \to TN$ whose kernel at every point is a Legendrian plane. The continuous map is just the retraction of $U_t$ onto $f_t(\Lambda)$, the bundle map is defined for points on $f_t(\Lambda)$ be being the identity map on $df_t(T\Lambda ) \sse TY$ to $T\Lambda \sse N$, by sending $P_t(\{0\} \oplus \R)$ isomorphically to $T(0,1) \sse TN$, and by letting $P_t(T^*\Lambda \oplus \{0\})$ be the kernel. We then extend $G_t$ to all of $U_t$ by parallel translation via some chosen contact connection.

We now cite \cite[Theorem 13.4.1]{EM} holonomic approximation for microflexible relations -- to construct a genuine Legendrian submersion. More specifically, we let $R_t \sse U_t$ be the set of all Reeb trajectories passing through $f_t(\Lambda)$. Then $f_t(\Lambda)$ is codimension $1$ in $R_t$, and therefore we can find a $C^0$--small isotopy of $f_t$ to $\wt f_t$, contained in $R_t$, and a Legendrian submersion defined near $\wt f_t$. Furthermore, the kernel of this submersion is $C^0$ close to the kernel of the formal submersion, $P_t(T^*\Lambda \oplus \{0\})$, which is transverse to $R_t$ everywhere and therefore in particular transverse to $\wt f_t(\Lambda)$. Here we are applying \cite[Theorem 13.4.1]{EM} in its relative parametric form: for $t \in \dd D^k$ where $f_t(\Lambda)$ is a Legendrian embedding the Legendrian submersion is equal to the front projection on the Weinstein Legendrian neighborhood. Additionally, the Reeb vector field remains transverse to $\wt f_t(\Lambda)$. This is because when using holonomic approximation the induces isotopies are graphical with respect to the positive codimension subset. This is explained in detail in \cite[8.3.1]{EM}.

So far, we have a family of smooth embeddings $\wt f_t: \Lambda \to Y$ so that near the image we have a Legendrian submersion to $\Lambda \x (0,1)$, so that the Legendrian fibers are transverse to $\wt f_t(\Lambda)$, as is the Reeb vector field. Let $\{(x_1, \ldots, x_n)\}$ be a local coordinate system for $\Lambda$, $\{(y_1, \ldots, y_n)\}$ be local coodinates on the Legendrian fibers, and let $z$ be a Reeb coordinate. Then together these form a local coordinate system for $Y$. Writing our contact form $\alpha$ in these coordinates, we see that $\alpha(\dd_z) = 1$ and $\alpha(\dd_{y_i}) = 0$ for all $i$. Thus $\alpha = dz - \sum g_i dx_i$, where $g_i: Y \to \R$ are some local functions. The contact condition is then equivalent to $\det\left[\frac{\dd g_i}{\dd y_j}\right] > 0$, and it follows using the implicit function theorem that we can find new coordinates $(\wt y_1, \ldots, \wt y_n)$ on the fibers so that $g_i = \wt y_i$. Thus, we have an explicit contactomorphism between a neighborhood of $f_t(\Lambda)$ and a neighborhood of some section of $\J^1\Lambda$.
\end{proof}

\section{Wrinkled Legendrians}\label{sec:wrinkled}

\subsection{Wrinkled embeddings}

We now discuss some of the model singularities needed to define wrinkled embeddings, and describe how they interact with Legendrian embeddings. Everything will be built from the model zig-zag. Define a plane curve $\psi: \R \to \R^2$ by 

$$\psi(u) = \left(\psi^x(u), \psi^z(u)\right) = \left(u^3 - u, \frac94u^5 - \frac52u^3 + \frac54u\right).$$

\begin{figure}
\includegraphics[scale = .3]{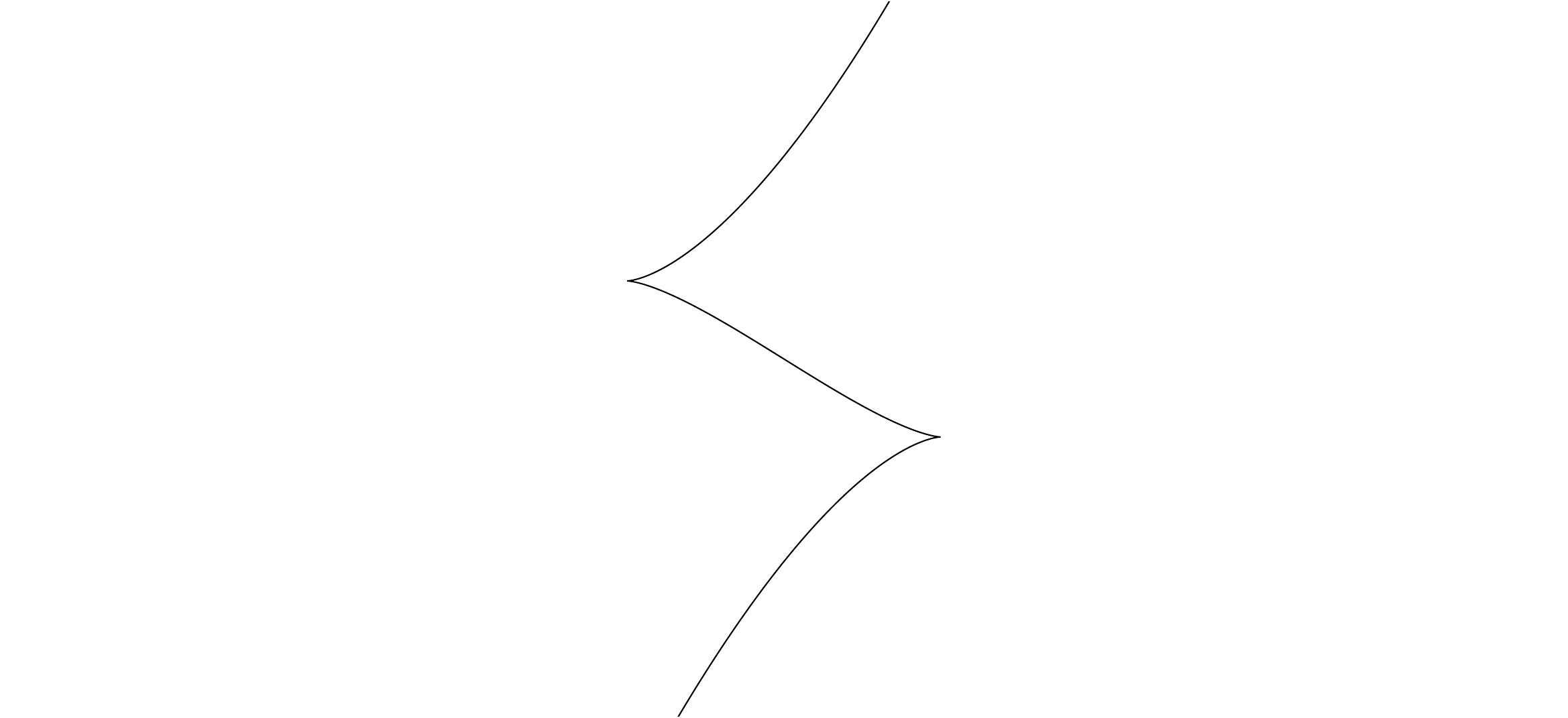} \caption{The curve $\psi$.}\label{fig:stab}
\end{figure}

The graph of $\psi$ is shown in Figure \ref{fig:stab}. Notice that $$\psi'(u) = \left(3\left(u^2 - \frac13\right), \frac{45}4\left(u^2 - \frac13\right)^2\right)$$
so in particular we see that $\psi$ is singular at $u = \pm \frac1{\sqrt3}$. By letting $\psi^y(u) = \frac{15}4(u^2-\frac13)$, we see that $\psi$ parametrizes the front projection of a smooth Legendrian curve in $\R^3_\std$. This will be relevant later.

We denote by $\psi_\delta$ a rescaling of $\psi$, defined by $$\psi_\delta(u) = \left(\delta^{\frac32}\psi^x\left(\frac u{\sqrt\delta}\right), \delta^{\frac52}\psi^z\left(\frac u{\sqrt\delta}\right)\right).$$

Notice that because of cancellation $\psi_\delta$ is defined for $\delta \leq 0$ as well. Explicitly: $$\psi_\delta(u) = \left(\psi_\delta^x(u), \psi_\delta^z(u)\right) = \left(u^3-\delta u, \frac94u^5 - \frac{5\delta}2u^3 + \frac{5\delta^2}4u\right).$$
 
Notice also that $\psi_\delta$ is a smooth graphical curve when $\delta < 0$. We will assume that each $\psi_\delta$ is a compactly supported curves -- this is false as written but it is easily accomplished with a cut-off function.

Consider the map $\R^n \to \R^{n+1}$ defined by \begin{equation}\label{eqn:ufst}(x_1, \ldots, x_{n-1}, u) \mapsto (x_1, \ldots, x_{n-1}, \psi^x_{x_{n-1}}(u), \psi^z_{x_{n-1}}(u)).\end{equation}

\begin{figure}
\includegraphics[scale = .7]{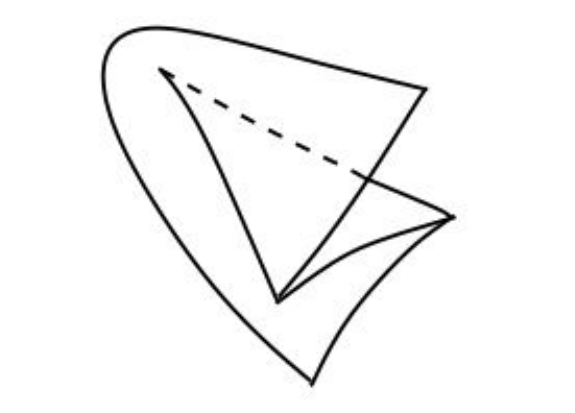} \caption{An unfurled swallowtail singularity.}\label{fig:3wrinkle}
\end{figure}

The image of this map is given in Figure \ref{fig:3wrinkle} in the $n=2$ case. The singular set is $\{3u^2=x_{n-1}\}$. When $x_{n-1}>0$ this singularity is diffeomorphic to the cusp singularity that is present for $\psi$, times $\R^{n-1}$. We call the singularity at $\{x_{n-1} = u = 0\}$ an \emph{unfurled swallowtail}. It is codimension $2$, and can be thought of as the singularity that occurs when $\psi$ is pulled tight, into a smooth graphical curve. Notice that an unfurled swallowtail is a topological embedding.

A \emph{wrinkle} is a map $w: \R^n \to \R^{n+1}$ given by $w(x, u) = (x, \psi_{1 - |x|^2}(u))$ (here $x \in \R^{n-1}$). Thus a wrinkle is singular on the sphere $\{|x|^2 + 3u^2 = 1\}$, there are cusp singularities on the upper and lower hemispheres and unfurled swallowtails along the equator $\{u=0\}$. An \emph{embryo} singularity is a singularity which is isolated in space and codimension $1$ in time, defined by the model $w_t(x, u) = (x, \psi_{t - |x|^2}(u))$ at $x=u=t=0$. Thus an embryo singularity is a singularity allowing wrinkles to appear or disappear in parametric families of maps (note that $w_t$ is a smooth map for $t < 0$). Alternatively, and embryo is just the singularity obtained when time-slices are tangent to unfurled swallowtail singularities in one larger dimension.

\begin{definition} \label{WE def}
Let $V$ and $W$ be manifolds, with $\dim W -1 = \dim V = n$. A \emph{wrinkled embedding} is a smooth map $f: V \to W$, which is a topological embedding, and which is singular on some finite collection of codimension $1$ spheres $S_j^{n-1} \sse V$, which bound disks $D^n_j \sse V$. Near each $S^{n-1}_j$ the map $f$ is required to be modeled on a wrinkle.\footnote{We do not require that $f$ is modeled on a wrinkle over all of $D^n_j$: wrinkles are allowed to be contained inside each other.}

In $k$--parametric families, wrinkled embeddings are allows to have embryo singularities, but no singularities of higher codimension.
\end{definition}

Notice that, although $\psi: \R \to \R^2$ is singular, the tangent plane to $\psi$ is well-defined everywhere. That is, if $G d\psi: \R \to \op{Gr}_{2,1}$ is the map taking a point in $\R$ to the image of $d\psi$ (a $1$--plane in $\R^2$), then even though $d\psi$ is $0$ at the singular points of $\psi$, $G d\psi$ extends to a smooth function over the entirety of $\R$. This is not true for arbitrary singularities! However, this is also true for the unfurled swallowtail, a fact which is easily checked (and we do so below). Thus, while a wrinkled embedding is singular at points, the map $G df: V \to \op{Gr}_{n+1, n}$ is smooth and defined everywhere.

These definitions are from \cite{wrinkled}. The main theorems there we will use are \cite[Theorem 2.5.1]{wrinkled} and its parametric version, \cite[Theorem 2.9.1]{EM}. We state it here for reference, specialized to the codimension $1$ case since this is all we will use.

\begin{thm}[\cite{wrinkled}]\label{thm:wrinkled}
Let $V$ and $W$ be manifolds, with $\dim W -1 = \dim V = n$. Let $f_t: V \to W$ be a parametric family of smooth embeddings, $t \in D^k$. Let $G_t^s: V \to \op{Gr}_n(W)$ be a smooth homotopy ($s \in [0,1]$) of maps covering $f_t$, so that $G_t^0 = G df_t$. Then there is a family of wrinkled embeddings $F_t^s: V \to W$, so that $F_t^0 = f_0$, and for all $s \in [0,1]$ $F_t^s$ is $C^0$--close to $f_t$, and $G dF_t^s$ is $C^0$--close to $G_t^s$ (with respect to a given local trivialization of $\op{Gr}_nW$). Furthermore, if on some closed set $A \sse V \x D^k$ the homotopy $G_t^s$ is constant in $s$, then on $A$ $F_t^s$ can be taken to be identically equal to $f_t$ as well. 
\end{thm}

As explained, this family of wrinkled embeddings has embryo singularities, but nothing more complicated. The set of $t \in D^k$ where embryos occur is a union of codimension $1$, smoothly embedded submanifolds of $D^k$. (In fact in \cite{wrinkled} it is shown that these submanifolds can be taken to be spheres, but we do not use this.)

\subsection{Wrinkled Legendrians}

We now explain how this relates to contact geometry. Consider the smooth manifold $W = \R \x \Lambda$. Contained in $\op{Gr}_nW$ is the set of $n$--planes which are transverse to the line field $T\R \x \{0\} \sse TW$, we call such planes \emph{non-vertical}. This subspace is canonically equivalent to $\J^1\Lambda$, as bundles over $\R \x \Lambda$: any non-vertical $n$--plane at $(z, x) \in \R \x \Lambda$ can be canonically written as the graph of a linear map $T\Lambda_x \to T\R_z$, and since $T\R_z$ is canonically isomorphic to $\R$ the space of such maps is equal to $T^*\Lambda_x$. This identification is natural with respect to contact geometry: if $V = \{(x,y,z) = (x, \sigma(x), h(x))\}$ for some sections $\sigma: \Lambda \to T^*\Lambda$ and $h: \Lambda \to \R$, then $V$ is a Legendrian submanifold if and only if $\sigma = Gd(\op{id} \x h)$ at every point $\x \in \Lambda$ (here $\op{id} \x h: \Lambda \to \Lambda \x \R$ is the obvious map).

In fact this correspondence goes beyond identification of Legendrian sections and holonomic sections: it holds for hypersurfaces in general. Suppose $f: V \to \Lambda \x \R$ is a map from any $n$--manifold $V$, so that $f$ is an immersion on some dense set $U \sse V$, and $G df: U \to \op{Gr}_n(\Lambda \x \R)$ is non-vertical at every point. If $Gdf$ extends smoothly to a function $G: V \to \op{Gr}_n(\Lambda \x \R)$ covering $f$, again non-vertical everywhere, then the smooth map $(f, G): V \to \J^1\Lambda$ is automatically a Legendrian immersion on $U$. Furthermore, as long as the map $(f,G)$ is a smooth immersion on all of $V$, it is automatically a Legendrian immersion, since being Legendrian is a closed condition and $U$ is dense.

The simplest example is the cusp singularity of a curve in $\R^2$: defining $f: \R \to \R^2$ by $f(t) = (x, z) = (t^2, t^3)$, we can compute $Gdf(t) = y = \frac32t$, simply by computing the slope $G = \frac{\dd f}{\dd z}\left(\frac{\dd f}{\dd x}\right)^{-1}$. The Legendrian condition forces this computation, and afterwards we see that $(f, G): \R \to \R^3$ is in fact an immersion. The computation above for $\psi^y(u)$ is similarly determined by $\psi^x(u)$ and $\psi^z(u)$, and we again have that the curve $u \mapsto (\psi^x(u), \psi^y(u), \psi^z(u))$ is smoothly embedded. We can either do this explicitly, or alternatively notice that the two singularities of the plane curve $\psi: \R \to \R^2$ are each diffeomorphic to the standard cusp (clearly the question of whether $(f, G)$ is an immersion for a given smooth singularity $f$ only depends on $f$ up to diffeomorphism).

Though the case of the cusp works out well, in general we are not so lucky: while the unfurled swallowtail and embryo singularities do lift canonically to smooth maps with target $\J^1\R^n$, the lifts are not immersions. Explicitly, given the unfurled swallowtail defined by Equation \ref{eqn:ufst}, we see that the lift is given by $y_i = 0$ for all $i \leq n-2$, $y_{n-1} = \frac52 u(x_{n-1} - u^2)$, and $y_n = \frac{15}4(u^2 - \frac{x_{n-1}}3)$. Over the set $\{x_{n-1} = u = 0\}$ (i.e. the unfurled swallowtail singularity) we see that $\dd_u$ is in the kernel of the smooth map $\R^n \to \R^{2n+1}$ defined by these coordinates. This describes the Legendrian lift of the singularity: it is smoothly embedded outside of the set $\{u=x_{n-1} = 0\}$, and on this set the differential has a rank $1$ kernel. Because they are built on the same singularity, this holds for wrinkles as well as embryo singularities.

We would like to use Theorem \ref{thm:wrinkled} to prove a statement about Legendrians, but because the singularities that appear in the statement do not arise as the front projection of any smooth Legendrian it does not directly apply. For the time being we will ``define away'' this problem.

\begin{definition} \label{def:WL}
Let $\Lambda$ be a smooth $n$-manifold and $(Y, \xi)$ a contact $(2n+1)$-manifold. A \emph{wrinkled Legendrian embedding} is a smooth map $f: \Lambda \to (Y, \xi)$, which is a topological embedding, satisfying the following properties. The image of $df$ is contained in $\xi$ everywhere, and $df$ is full rank outside a subset of codimension $2$. This singular set is required to be diffeomorphic to a disjoint union of $(n-2)$-spheres $\{S_j^{n-2}\}$, called \emph{Legendrian wrinkles}. Each $S^{n-2}_j$ is contained in a Darboux chart $U_j$, so that $\Lambda \cap U_j$ is diffeomorphic to $\R^n$, and the front projection $\pi_j \circ f: \Lambda \cap U_j \to \R^{n+1}$ of $f$ is a wrinkled embedding, smooth outside of a compact set. In particular, the front projection of each $S^{n-2}_j$ is the unfurled swallowtail singular set of a single wrinkle in the front.

For parametric families of wrinkled Legendrians we also allow \emph{Legendrian embryos}; Legendrian lifts from the front projection of embryo singularities.
\end{definition}

A wrinkle Legendrian is therefore a smooth Legendrian embedding outside a set of codimension $2$, however it is permitted to contain the singularity defined as the Legendrian lift of the unfurled swallowtail. Our definition is slightly stronger than this: we also require a global trivialization of each Legendrian wrinkle given by the Darboux charts $U_j$. We emphasize that $\{U_j\}$ is considered part of the data of a wrinkled Legendrian: for a given map $f:\Lambda \to (Y, \xi)$, different choices of $\{U_j\}$ are considered to be different as wrinkled Legendrians. However notice there is no requirement for these Darboux charts to be disjoint, and in fact we often take them to be equal when multiple Legendrian wrinkles are contained in a single Darboux chart.

When a Legendrian wrinkle is born we add a new $U_j$ to the collection of Darboux charts which contains the Legendrian embryo, and it is required to contain the created wrinkle throughout its entire ``lifetime''. To topologize the space of wrinkled Legendrians we use the $C^\infty$ topology on the space of maps, together with independent $C^\infty$ topologies for the Darboux charts $U_j$, with the proviso that a $U_j$ is allowed continuously appear or disappear at Legendrian embryos.

Combining Proposition \ref{prop:graphical} with Theorem \ref{thm:wrinkled} then immediately gives us a theroem about wrinkled Legendrians. However, to state it sensibly, we need to define a map from wrinkled Legendrian embeddings to formal Legendrian embeddings. It is not immediately obvious how to do this, since for a wrinkled Legendrian $f: \Lambda \to (Y, \xi)$ the differential $df$ is not full rank everywhere. This is where we will use the data of the Darboux charts $U_j$ in the definition of a wrinkled Legendrian.

Outside of our Darboux charts $f$ is already Legendrian, therefore we just define the formal Legendrian corresponding to a wrinkled Legendrian on a given Darboux chart. This is essentially the same as \emph{regularization} from \cite{wrinkled}. On a given chart $U_j$, we first perturb our map $f: \R^n \to U_j$ to a smooth embedding by fixing the $y$--coordinates and perturbing its front projection to a smooth embedding $\R^n \to \R^{n+1}$. There is a natural way to do this, just rounding out all the singularities. The resulting smooth embedding $\wt f$ is not Legendrian, but we claim that $d\wt f$ is homotopic to a Legendrian bundle map in a canonical way. This homotopy again fixes the $\dd_{y_i}$ coordinates of the output of $d\wt f$, and rotates the map in the $x_n$-$z$ plane until the map is everywhere transverse to $\dd_z$ (in the coordinates given by Equation \ref{eqn:ufst}). Once the map is transverse to $\dd_{y_i}$ and $\dd_z$, we note that the space of all maps transverse to this $(n+1)$--plane field is contractible (being homeomorphic to the space of $(n+1) \x n$ matrices), and the set of bundle maps inside this space with Legendrian image is also contractible (being homeomorphic to the space of symmetric $n \x n$ matrices). Therefore we can find a homotopy through bundle maps transverse to $\dd_{y_i}$ and $\dd_z$ to a bundle map with Legendrian image, and this homotopy is canonical (up to further homotopy).

Altogether this defines a map from wrinkled Legendrian embeddings to formal Legendrian embeddings. While the explicit map depends on many choices, the map is canonical up to homotopy, and in particular can be made continuous in any $k$--parametric family.

\begin{prop}\label{prop:wrink Leg}
Let $(f_t, F_{s,t})$ be a parametric family of formal Legendrian embeddings $\Lambda \to (Y, \xi)$, $t \in D^k$. Then the family $(f_t, F_{s,t})$ is homotopic through formal Legendrian embeddings to a family $\ol f_t: \Lambda \to (Y, \xi)$ of wrinkled Legendrian embeddings. If $(f_t, F_{s,t})$ is already a wrinkled Legendrian embedding on a closed subset $A \sse \Lambda \x D^k$, then we can take $\ol f_t = f_t$ on this set.
\end{prop}

The proof is just an immediate combination of Proposition \ref{prop:graphical} with Theorem \ref{thm:wrinkled}. We also note that, since Proposition \ref{prop:graphical} and Theorem \ref{thm:wrinkled} are both $C^0$--dense, Proposition \ref{prop:wrink Leg} is $C^0$--dense as well.

We will mostly be interested in applications to the case when $n \geq 2$, because in these dimensions there are additional surgery techniques we can use. When $n=1$: notice that the definition of a wrinkled Legendrian embedding becomes somewhat strange: a single wrinkled Legendrian curve is just a smooth Legendrian curve, since the allowed singularities are in codimension $2$. However, a difference arises when considering families of wrinkled Legendrian curves, which are allowed to have embryo singularities. Passing through an embryo singularity is locally contactomorphic to the family $(x, z) = \psi_t(u)$; that is, it takes a curve with smooth front projection and introduces a zig-zag (or deletes a zig-zag that already exists). In the literature this operation is called the \emph{stabilization} of a Legendrian knot, and in this way Proposition \ref{prop:wrink Leg} recovers a theorem of Fuchs and Tabachnikov.

\begin{thm}[\cite{FT}]
Let $f_t: S^1 \to (Y, \xi)$ be a family of Legendrian embeddings, with $t \in S^{k-1}$, so that the family $\{f_t\}$ extends to $D^k$ among formal Legendrian embeddings. Then, after stabilizing the Legendrians $f_t$ sufficiently many times (with both orientations), the family extends to $D^k$ among genuine Legendrian embeddings.
\end{thm}

When $k=1$, this briefly states that formally isotopic Legendrian knots become Legendrian isotopic after some number of stabilizations.\footnote{It is furthermore true that smoothly isotopic Legendrian knots become formally isotopic after sufficiently many stabilizations. This fact does not generalize for $k>1$, however.} It was shown by Etnyre and Honda \cite{EH} that in fact the number of stabilizations required is unbounded when considering all pairs of formally equivalent Legendrians. In a sense, this is the principal difference between knots and the $n\geq 2$ case.

\begin{remark}
When $k=0$, it may appear that the result claims that any formal Legendrian embedding is isotopic to a genuine Legendrian embedding, but in fact this result is false \cite{Benn}. The discrepancy comes from our identification of wrinkled Legendrians with smooth Legendrians: while the underlying embedding of a wrinkled Legendrian is a smooth Legendrian embedding, the formal isotopy class represented by a wrinkled Legendrian depends on both the underlying smooth map as well as the Darboux charts $\{U_j\}$ covering the wrinkles. Indeed, the definition of a formal Legendrian associated to a wrinkled Legendrian is exactly a formal destabilization of the bundle map near the wrinkle. Thus the $k=0$ case states that any formal Legendrian knot is a formal destabilization of a genuine Legendrian knot.
\end{remark}

\section{Loose Legendrians}\label{sec:loose}

\subsection{Surgery of singularities}

In order to make Proposition \ref{prop:wrink Leg} more useful, we would like a way to remove the singularities from the family. In general there is no way to do this consistently. However, if we are working with something which is already singular, we can use surgery of singularities in order to get one large connected singular set which never needs to change throughout the family. Similarly, if we have something which is smooth but ``looks like'' a resolved singularity, we will be able to use it as a substitute for ever needing a singular set. This will lead us to the definition of loose Legendrians.

For any $\delta > 0$, we fix the notation $m_\delta: \R \to (0, \infty)$ to be a smoothing of the function $x \mapsto \min(x, \delta)$. Specifically, we ask that $m_\delta(x) = x$ for all $x > 2\delta$, $m_\delta(x) = \delta$ for all $x \leq \frac12\delta$, $\delta \leq m_\delta(x) \leq x$ when $\frac12\delta < x < 2\delta$, and $m_\delta'(x) \leq 1$ everywhere. 

\begin{definition}
Let $f: \Lambda \to (Y, \xi)$ be a wrinkled Legendrian. A \emph{marking} for $f$ is a compact codimension $1$ embedded submanifold $\Phi \sse \Lambda$, so that the boundary of $\Phi$ is a disjoint union of spheres which are mapped via $f$ to a subset of the Legendrian wrinkles. We require that in the local model $(\pi_j \circ f)(\Lambda \cap U_j) \cong w(\R^n)$, $\Phi$ is given as $\Phi = \{u=0, |x| \geq 1\}$, and that the interior of $\Phi$ is disjoint from the singular set of $f$.

In parametric families $f_t: \Lambda \to (Y, \xi)$, we require the family $\Phi_t \sse \Lambda$ to be smoothly varying in $t$ whenever we are disjoint from the set of $t \in D^k$ corresponding to embryo singularities. At an embryo singularity modeled in the front projection $(\pi_j \circ f_s)(\Lambda \cap U_j) \cong w_s(\R^n)$, we require $\Phi_s = \{u=0, |x|^2 \geq s\}$. Here $s \in (-\e, \e)$ is a chosen local coordinate in $D^k$ which is transverse to the embryo singular set, which is a smooth codimension $1$ submanifold (indeed \cite{wrinkled} guarantees it can be taken to be a sphere though we do not use this anywhere). In particular, the topology of $\Phi_t$ is allowed to change at embryo singularities, with $\Phi_s$ for $s>0$ being abstractly diffeomorphic to $\Phi_s$ for $s<0$ with an open disk removed, so that the new boundary component is exactly the new singular sphere of $f_s$. We also allow the possiblility for $\Phi_t$ to be disjoint from the embryo in which case it remains disjoint from the new singular sphere.
\end{definition}

\begin{figure}
\includegraphics[scale=.5]{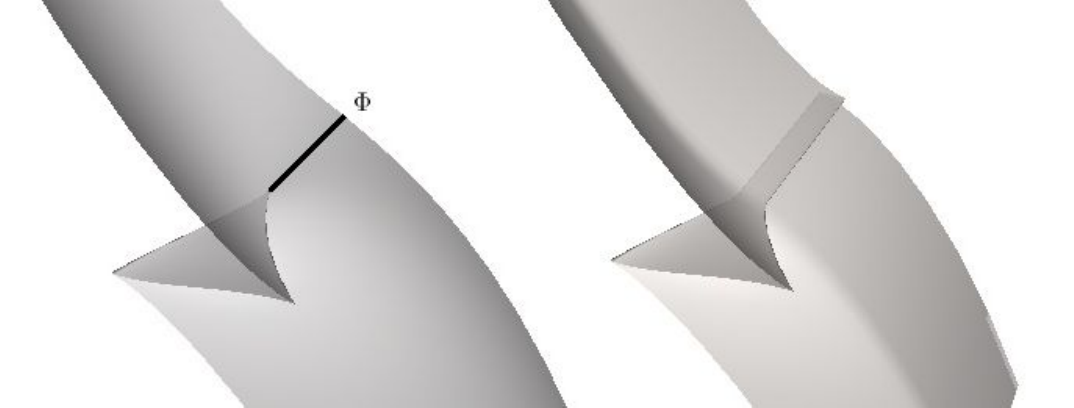} \caption{A local model describing how a marking resolves Legendrian wrinkles.} \label{fig:resolve}
\end{figure}

The idea of a marking is that it is a pattern that gives us a canonical way to desingularize wrinkles. Intuitively, the singular set of $f$ looks like zig-zags in the front which are pulled tight to become graphical; the zig-zags themselves are smooth but the location where they are pulled tight is singular. Instead, we can just let the zig-zags persist as tiny zig-zags throughout $\Lambda$. The marking $\Phi$ is just some formal data that tells us where to put the tiny zig-zags in a consistent way. See Figure \ref{fig:resolve}.

To be more formal, let $u \in (-1, 1)$ be a coordinate in $\Lambda$ transverse to $\Phi$ and let $\wh \Phi = \Phi \cup_\dd [0, \e)$, so that $\wh \Phi \x (-1, 1)$ is an open neighborhood of $\Phi$ in $\Lambda$. Let $r: \wh\Phi \to \R$ be the distance to the boundary, defined so that $r$ is negative exactly on $\wh \Phi \sm \Phi$. Then by the Legendrian tubular neighborhood theorem and our assumption on $\Phi$ there is a neighborhood $U_\Phi \sse Y$ of $f(\Phi)$ so that $U_\Phi$ is contactomorphic to $\J^1(\wh \Phi \x (-\e, \e))$ and the front of $f(\Lambda) \cap U_\Phi$ is given by $(x, z) = \psi_{-r}(u)$. Here $x \in (-\e, \e)$ and the front is the identity map on the $\wh \Phi$ component.

Then, we simply replace $f$ with $\wt f$, defined by $(x, z) = \psi_{m_\delta(-r)}(u)$, where $\delta$ is chosen to be much smaller than $\e$. The resulting front has only cusp singularities, therefore $\wt f$ is smooth everywhere on a neighborhood of $\Phi$, and equal to $f$ outside of this neighborhood. It is also clear we can do this parametrically, since at embryos the model $\Phi_t$ is set up to be consistent with the model embryo singularity. Thus we have shown:

\begin{lemma}\label{lem:resolve}
Let $f_t: \Lambda \to (Y, \xi)$ be a family of wrinkled Legendrians with $t \in D^k$, and let $\Phi_t \sse \Lambda$ be a family of markings for $f_t$. Then $f_t$ is $C^0$-close to a wrinkled Legendrian family $\wt f_t: \Lambda \to (Y, \xi)$, so that $\wt f_t$ is smooth on a neighborhood of $Phi_t$ and equal to $f_t$ outside that neighborhood.
\end{lemma}

So far, if we want to construct a smooth family of Legendrians, it only remains to find a marking. In fact this is easy, but a difficulty arises in the relative parametric case. Though the topology of $\Phi_t$ is allowed to change, it changes only by removing open disks, and therefore if $\Phi_t$ is nonempty for a single $t$ it will be non-empty for all $t$. In particular, if $f_t$ is smooth for $t \in \dd D^k$ but singular on the interior of $D^k$, then any possible $\Phi_t$ we might choose to resolve our singularities will consist of closed manifolds for $t \in \dd D^k$. Of course, even if $f$ is smooth, $\wt f$ will not generally be equal to $f$ unless $\Phi = \es$. This will motivate our definition of loose Legendrians: intuitively they are smooth Legendrians which look like resolutions of some wrinkled Legendrian along some marking.

\subsection{Loose Legendrians}

Let $C \sse \R^3$ be the cube of side length $1$, and let $\Lambda_0 \sse C$ be a properly embedded Legendrian arc whose front is a zig-zag and which is equal to the set $\{y=z=0\}$ near the boundary. Let $V_\rho$ be the open set $\{|q| < \rho, |p| < \rho\} \sse T^*\R^{n-1} = \R^{n-1} \x \R^{n-1}$, and let $Z_\rho \sse V_\rho$ be the intersection of $V_\rho$ with the zero section $\{p = 0\} \sse T^*\R^{n-1}$. Letting $\lambda = \sum p_idq_i$ and defining $\alpha_\std = \alpha_0 - \lambda$ on $\R^3 \x T^*\R^{n-1}$, we see that $\alpha_\std$ defines the standard contact structure on $\R^3 \x T^*\R^{n-1} = \R^{2n+1}_\std$. In these coordinates we see $C \x V_\rho$ as an open set in $\R^{2n+1}_\std$, and $\Lambda_0 \x Z_\rho$ is a Legendrian submanifold in this set.

\begin{definition}
Let $\rho > 1$. Then any Legendrian $\Lambda_\ell$ is a contact manifold $U$ which is contactomorphic as a pair to $(C \x V_\rho, \Lambda_0 \x Z_\rho)$ is called a \emph{loose chart}.

If $\Lambda \sse (Y, \xi)$ be a Legendrian submanifold. Then we say that $\Lambda$ is \emph{loose} if there is an open set $U \sse Y$ so that $(U, U \cap \Lambda)$ is a loose chart.
\end{definition}

\begin{prop}\label{prop: two loose}
For any $\ol \rho > 1$, an arbitrary loose chart contains another loose chart of size parameter $\ol\rho$. In particular a loose chart contains two disjointly embeddded loose charts, and therefore a loose chart contains infinitely many disjointly embedded loose charts.
\end{prop}

\begin{proof}
For $h \in (0,1)$ let $h\Lambda_0$ be the Legendrian curve obtained by scaling $\Lambda_0$ via the contactomorphism $(x,y,z) \mapsto (hx, hy, h^2z)$. Notice that $h \Lambda_0 = \Lambda_0$ near the boundary, and that they are Legendrian isotopic rel boundary.

Inside $C \x V_\rho$, define the Legendrian $\Lambda$ to be the Legendrian whose front is given as $\{(q,x,z); (x,z) \in m_\delta(|q| +1 - \rho)\Lambda$. Notice that $\Lambda \sse C \x V_\rho$ since $|p| = m'_\delta(|q| +1 - \rho) \leq 1$. Then $\Lambda$ is equal to $\Lambda_0 \x Z_\rho$ near their boundaries, and they are Legendrian isotopic rel boundary since their fronts are smoothly ambient isotopic rel boundary. Thus there is a compactly supported contact isotopy taking $\Lambda_0 \x Z_\rho$ to $\Lambda$, and therefore it suffices to find two disjoint loose charts in $(C \x V_\rho, \Lambda)$. Let $\wt \rho = \rho - 1 - \frac12 \delta$, and restrict to the subset $|q| < \wt \rho$. On this set $\Lambda = \delta\Lambda_0 \x Z_{\wt \rho}$. Scaling via the contactomorphism $(q, p, x, y, z) \mapsto (\frac q\delta, \frac p \delta, \frac x \delta, \frac y\delta, \frac z{\delta^2})$, we see that this set is contactomorphic to $\Lambda_0 \x Z_{\frac{\wt\rho}\delta} \sse C \x V_{\frac{\wt\rho}\delta}$. From here, we simply note that by choosing $\delta$ to be small we can make $\frac{\wt\rho}{\delta}$ as large as we please.
\end{proof}

Define an \emph{inside-out wrinkle} to be the map $\ol w: \R^n \to \R^{n+1}$ defined by $\ol w(x, u) = (x, \psi_{|x|^2-1}(u))$. Then $\ol w$ is singular on the hyperbola $\{|x^2| - 3u^2 = 1\}$, which has unfurled swallowtails on the subset $\{u=1-|x|^2=0\}$ and cusps elsewhere. In particular it is not smooth even outside of a compact set. The wrinkled Legendrian $f: \R^n \to \J^1\R^n$ whose front projection is $\ol w$ is smooth outside of a compact set (since the front has only cusp singularities), but it is not standard.

However, let $\Phi = \{u=0, |x| \leq 1\}$. Then $\Phi$ is a marking for $f$. Furthermore, if we resolve $f$ along $\Phi$ to obtain $\wt f$, then $\wt f(\R^n) \cap B^{2n+1}(\rho)$ is contactomorphic to a loose chart for any $\rho > 1$: for any fixed $x_0 \in \R^{n-1}$, $\wt f(\R^n) \cap \{x=x_0\}$ is a zig-zag, therefore we can isotope $\wt f(\R^n) \cap B^{2n+1}(\rho)$ rel boundary to be a trivial family of standard zig-zags.

\begin{proof}[Proof of Theorem \ref{thm:main}:]
First, we apply Proposition \ref{prop:wrink Leg} to replace our family of formal Legendrian embeddings by a family of wrinkled Legendrian embeddings, which is equal to the original genuine Legendrian family on $(\Lambda \x \dd D^k) \cup (U \x D^k)$; we continue to denote this family by $f_t: \Lambda \to Y$. The set of $t \in D^k$ where embryo singularities appear is a union of smoothly embedded codimension $1$ submanifolds; let $K \in \Z_+$ be the number of these manifolds. By assumption $f_t^{-1}(U) = D^n$ is a fixed disk and $f_t|_{D^n}$ is constant in $t$ and $(U, f_t(D^n))$ is a fixed loose chart. Choose disjoint open sets $U_i \sse U$, $i=1, \ldots, K$, so that each $f_t(\Lambda) \cap U_i$ is a loose chart, which exist following Proposition \ref{prop: two loose}.

We define a new family $g_t: \Lambda \to Y$ of wrinkled Legendrians which is equal to $f_t$ outside of $\bigcup_{i=1}^K U_i$, and on each $f_t^{-1}(U_i)$ we replace the loose chart with an inside-out wrinkle with the same boundary conditions; again constant in $t$. For each $i$, we define a marking $\Phi^i_t$ of $g_t$ so that:

\begin{itemize}
\item near $\dd D^k$, $\Phi^i_t$ is a disk contained in $g_t^{-1}(U_i)$, exactly as defined above for the model inside-out wrinkle,
\item for all $t \in D^k$, $\Phi_t^i$ is either diffeomorphic to a disk $D^{n-1}$ or a cylinder $S^{n-1} \x [0,1]$,
\item for the $i$th connected component of the embryo set, these embryos are all contained in $\Phi_t^i$, and
\item other than the component in $U_i$, the boundary of $\Phi_t^i$ is exactly the sphere of Legendrian wrinkles created by this embryo set.
\end{itemize}

This can obviously be done: inside $\Lambda$ the singular set is a union of contractible embedded codimension $2$ spheres; any such family of spheres can always be realizes as the boundary of a family of embedded disks (or a family of embedded cylinders whose other boundary is another fixed contractible sphere). By definition of a marking we are also requiring $\Phi_t^i$ to be disjoint from all other components of the singular set, which is easy to arrange. We do not require the $\Phi_t^i$ to be disjoint for different $i$, though this can also be done in general.

We now just apply Lemma \ref{lem:resolve} for each $\Phi_t^i$, one at a time. The resulting family $\wt g_t: \Lambda \to (Y, \xi)$ is a family of smooth Legendrians. Furthermore, for $t \in \dd D^k$ $\wt g_t$ is isotopic to $f_t$, via an isotopy supported in $\cup_i U_i$. Thus we can add a collar isotopy (constant in $t$) between $\{f_t\}_{t \in \dd D^k}$ and $\{\wt g_t\}_{t \in D^k}$, the union of this collar and $\{\wt g_t\}_{t \in D^k}$ is a family of genuine Legendrians extending the family $\{f_t\}_{t \in \dd D^k}$ over the disk.

To conclude the proof of Theorem \ref{thm:main}, it remains to show that this extension is formally Legendrian isotopic through $D^k$ families rel $\dd D^k$ to the original $f_t$. This is not completely obvious, since we have not defined a formal Legendrian embedding corresponding to $g_t$. There are two issues: there is no defined formal Legendrian embedding corresponding to an inside-out wrinkle, and we have not understood how resolving singularities via a marking $Phi$ affects the formal Legendrian isotopy class. In a sense these two issues cancel each other out, since we can compare $f_t$ to $\wt g_t$ directly. Suppose $\Phi^i_t$ is diffeomorphic to a cylinder for a fixed $t$ and consider a neighborhood in $Y$ which is the union of three open sets: the neighborhood of $\Phi_t$ in the proof of Lemma \ref{lem:resolve}, the loose chart $U_i$ containing an inside-out wrinkle, and the Darboux chart containing the wrinkle on the other boundary component of $\Phi^i_t$. The front of $g_t$ in this neighborhood has the model singularities at the wrinkle and the inside-out wrinkle, but it is smooth elsewhere, since $g_t$ is equal to the zero section near the interior of $\Phi^i_t$.

Now, if we resolve $g_t$ along the small disk marking of $U_i$, we get the original wrinkled Legendrian $f_t$. In our front projection, we have a $D^{n-1}$ family of zig-zags in the loose chart, and a disjoint wrinkle. This is formally Legendrian isotopic to a loose chart, since by definition a wrinkle is formally Legendrian isotopic to the zero section. On the other hand if we resolve $g_t$ along $\Phi^i_t$ we obtain $\wt g_t$. In our front projection here, we still just have a $D^{n-1}$ family of zig-zags which is and a smooth front everywhere else; therefore this is formally Legendrian isotopic to a loose chart as well. Therefore $f_t$ is formally Legendrian isotopic to $\wt g_t$. 

Since this formal Legendrian isotopy between $f_t$ and $\wt g_t$ is canonical up to contractible choices, it is continuous in $t$. We also observe that resolving $g_t$ along $\Phi_t$ either when $\Phi^i_t$ contains an embryo or when $\Phi_t^i$ is a disk containing no singularities besides the inside-out wrinkle, the result is the same: a front which is smooth except for a disk of zig-zags.
\end{proof}

\section{Conclusion}\label{sec:conclusion}

We now state some basic applications of Theorem \ref{thm:main}, as well as some related results. The first is Theorem \ref{thm:main iso}.

\begin{proof}[Proof of Theorem \ref{thm:main iso}:]
We would like to apply Theorem \ref{thm:main} in the case $k=1$, but the remaning difficulty is that Theorem \ref{thm:main iso} does not assume that the loose charts for $f_0$ and $f_1$ are identical. It is a basic fact of contact topology that any two Darboux charts are ambient contact isotopic, even as parametrized open balls. Consider the new isotopy $\wt f_t = \p_t \circ f_t$, then there is a fixed Darboux chart $U \sse Y$ so that $\wt f_0^{-1}(U)=\wt f_1^{-1}(U) = D^n \sse \Lambda$, $\wt f_0 = \wt f_1$ on $D^n$, and the pair $(U, \wt f_i(D^n))$ is a loose chart for $i=0,1$. Next, smoothly isotope the isotopy $\wt f_t$ rel $\dd D^1$ to a new smooth isotopy $g_t$, so that $g_t^{-1}(U) = D^n$ and $g_t|_{D^n} = \wt f_0$ for all $t \in D^1$. Since the space of formal Legendrian  embeddings is a Serre fibration over the set of smooth embeddings, $g_t$ can be realized as a formal Legendrian embedding $(g_t, G_{s,t})$, and this formal Legendrian isotopy satisfies the hypothesis of Theorem \ref{thm:main}.
\end{proof}

We also have the corollary to Theorem \ref{thm:main} for the case $k=0$.

\begin{coro}\label{cor:exist}
Suppose $n \geq 2$ and let $(f, F_s)$ be a formal Legendrian embedding $\Lambda \to (Y, \xi)$. Then $(f, F_s)$ is formally Legendrian isotopic to a genuine Legendrian $\wt f: \Lambda \to (Y, \xi)$ (which is unique up to Legendrian isotopy by Theorem \ref{thm:main iso}).
\end{coro}

One proof of this fact can be done as follows. First, prove that any formal Legendrian embedding is isotopic to some Legendrian embedding. Then, find an explicit loose Legendrian in $\J^1(\R^n)$, which can be taken to the zero section $\{z = y = 0\}$ via an compactly supported formal Legendrian isotopy. Since this zero section is the standard local model of any Legendrian, we can there therefore replace any small open disk in a given Legendrian with a formally isotopic disk which is loose.

The existence of such a loose Legendrian plane can be proved directly quite easily, using for example stabilization near a cusp as in \cite{Stein} \cite{EES}. To prove that any formal Legendrian embedding is formally isotopic to a genuine Legendrian embedding can be proven by using an $h$--principle for Legendrian immersions, then noting that a generic immersion is embedded and classifying formal Legendrian embeddings in a formal Legendrian regular homotopy class.

The proof we present here is essentially the $k=0$ case of Theorem \ref{thm:main}. It is fairly distinct from the proof above, in that it does not use the $h$--principle for Legendrian immersions, and it requires no classification of formal Legendrian isotopy classes.

\begin{proof}
First by smooth isotopy we can assume that there is a set $U$ so that $f$ is Legendrian on $f^{-1}(U)$ and that $(U, f(\Lambda) \cap U)$ is a loose chart. From there, we use Proposition \ref{prop:wrink Leg} (relative to $U$) to replace $f$ by a wrinkled Legendrian embedding. The singular set is made up of $K$ disjoint copies of $S^{n-2}$ in $\Lambda$. Find $K$ loose charts inside $U$, and for each of these loose charts replace the Legendrian with an inside-out wrinkle. For each $i=1,\ldots, K$, find a marking $\Phi^i$ which is a cylinder joining one of the inside out wrinkles to one of the wrinkles outside $U$. Then apply Lemma \ref{lem:resolve}.
\end{proof}  

Demonstrating the existence of non-loose Legendrians requires some technology, but it is a simple result with the correct tools. Using Legendrian contact homology, we can easily see that the $LCH$ algebra of a loose Legendrian vanishes, since after isotopy we can find a Reeb chord bounding a single pseudo-holomorphic half-disk which has smaller action than all other chords, see \cite{EES}. This shows that any Legendrian with non-vanishing $LCH$ is non-loose. In particular any Legendrian with an augmentation is non-loose, so any Legendrian which is the convex boundary of an exact filling is non-loose as well \cite{E}. It is also easy to see that a loose Legendrian in any co-sphere bundle is never the singular support of a constructable sheaf: after Legendrian isotopy we can always arrange that the front contains a zig-zag

We emphasize that Theorem \ref{thm:main iso} is a staement about parametrized Legendrians. In particular, if $\p: \Lambda \to \Lambda$ is a diffeomorphism which is homotopic to the identity, then for any Legendrian embedding $f: \Lambda \to (Y, \xi)$, we can easily see that $f$ is formally Legendrian isotopic to $f\circ \p$ (since we are in high dimensions we can use the Whitney trick to find a smooth isotopy, then using obstruction theory we see this can be made into a formal Legendrian isotopy since $\p^*: H^*(\Lambda) \to H^*(\Lambda)$ is the identity). Thus is follows that, if $f$ is loose, then $f \circ \p$ is loose as well (since looseness is a property of unparametrized submanifolds), and therefore $f$ and $f \circ \p$ are Legendrian isotopic. This is known to be false in general: the first such example is due to Abouzaid \cite{Ab}, where he shows in particular that the standard Legendrian sphere $S^{4k+1} \sse S^{8m +3}_\std$, reparametrized by a diffeomorphism which defines an exotic sphere which does not bound a parallelizable manifold, is not isotopic to the same sphere with the standard parametrization.

We emphasize that the essential property of a loose chart is the condition $\rho>1$: in fact for every Legendrian submanifold $\Lambda$ of dimension $n \geq 2$, we can find a Darboux chart $(U, U \cap \Lambda) \cong (C \x V_\rho, \Lambda_0 \x Z_\rho)$ for any $\rho < 1$. To show that such a Darboux chart exists for some $\rho <1$ we can appeal to an $h$--principle of Gromov's \cite{PDR} about open contact submanifolds of positive codimension (the manifold $C$ in this case). To see the result for all $\rho < 1$ it suffices to do it for the zero section, which can be done via an explicit isotopy.

Since these size issues are fairly subtle one might imagine it would be very difficult in practice to determine whether a given Legendrian is loose or not. However, Legendrians are typically presented via a front projection, and in this case there is an easy criterion.\footnote{This proposition was also proved in \cite{affine}.}

\begin{prop}
Let $\Lambda \sse \J^1(Q)$ be a Legendrian submanifold, and let $\pi(\Lambda)$ be its front projection. Let $D^2 \sse \R \x Q$ be a disk which is tangent to the vertical direction, and suppose that $D^2$ is transverse to $\pi(\Lambda)$, meaning that it is transverse to the smooth stratum of $\pi(\Lambda)$, transverse to the sets of cusps in $\pi(\Lambda)$, and disjoint from all singularities of codimension $2$ or larger. Suppose that $D^2 \cap \pi(\Lambda)$ is equal (as a curve in $D^2$) to a zig-zag which is equal to a horizontal line near $\dd D^2$. Then $\Lambda$ is loose.
\end{prop}

\begin{proof}
Let $W = D^2 \x B^{n-1}$ be an open neighborhood of $D^2$, with $B^{n-1}$ chosen small enough so that $W \cap \pi(\Lambda)$ is diffeomorphic to $(D^2 \cap \pi(\Lambda)) \x B^{n-1}$. Let $U \sse \J^1(Q)$ be the set of points projecting to $W$. Then $U$ is contactomorphic to $C \x T^*B^{n-1}$, and under this contactomorphism $\Lambda \cap U$ is equal to $\Lambda_0 \x \{p=0\}$. From here, we take the symplectomorphism $T^*B^{n-1} = \{|q| < \e\} \sse T^*\R^{n-1}$ given by $(q, p) \mapsto (\frac2\e q, \frac\e 2 p)$. This map sends $T^*B^{n-1}$ to $\{|q| < 2\}$, furthermore it preserves the set $\{p=0\}$, and it preserves $\lambda$, therefore it defines a contactomorphism from $C \x T^*B^{n-1}$ to $C \x \{|q|<2\}$ by acting by the identity on $C$, sending $\Lambda_0 \x \{p=0\}$ to $\Lambda_0 \x \{p=0\}$. This latter set contains a loose chart.
\end{proof}

In fact the hypotheses of the proposition are stronger than necessary: any $D^2$ suffices whether it is tangent to the vertical or not, and as long as $D^2 \cap \Lambda$ is diffeomorphic as a curve in $D^2$ to some zig-zag the theorem holds. This follows because ambient isotopies of $\R \x Q$ give fronts which define isotopic Legendrians, as long as those fronts are never tangent to the vertical. Note however that there are some subtleties in other cases, for example when considering fronts in $\J^1(Q, S^1)$.

\appendix

\section{Invariants of formal Legendrian isotopy classes} \label{appx}

In order for the main result of this paper to useful in practice, we would like to have an explicit way to tell when two knots are formally isotopic. We do this for the case of Legendrians in $\R^{2n+1}_\std$, similar computations can be done in other manifolds but the algebraic topology becomes more difficult. Up to some $\Z_2$ indeterminacy, it turns out that formal isotopy classes are in correspondence with two invariants we understand well: the classical invariants $tb$ and $r$. Some of the details in calculation are left to the reader, they can also be found in ~\cite{EES}.

\begin{definition}
Let $(f, F_s)$ be a formal Legendrian embedding in $(Y,\xi)$. $F_1$ is a bundle map $T\Lambda \to f^*\xi$, so every fiber has Lagrangian image. The homotopy class of this map in the space of Lagrangian bundle monomorphisms $T\Lambda \to f^*\xi$ is called the \emph{rotation class} of $(f, F_s)$. We denote this class $r(f, F_s)$.
\end{definition}

Immersed Legendrians satisfy an $h$--principle ~\cite{PDR}, and the rotation class classifies them up to regular Legendrian homotopy. If we have two formal Legendrian embeddings which are smoothly isotopic, we can compare their rotation classes as follows. $F_1$ defines an isomorphism $f^*\xi \cong T\Lambda \otimes \C$, therefore two formal Legendrians together define an element in $\op{Aut}_{\C}(T\Lambda \otimes \C)$. Two Legendrians have the same rotation class if and only if this difference element is in the component of the identity. If $f^*\xi$ is trivial (which is always the case if $\xi$ is a trivial bundle on $Y$) then $\op{Aut}_{\C}(f^*\xi) \cong \op{Map}(\Lambda, U_n)$, thus the difference class $r(f, F_s) - r(\wt f, \wt F_s)$ is an element of $K^1(\Lambda)$ in this case.

\begin{definition}
Suppose $n$ is odd, and let $(f, F_s)$ be a formal Legendrian knot in $(Y, \xi)$. Assume $\Lambda$ is orientable, and that $f(\Lambda)$ is nulhomologous and coorientable in $(Y, \xi)$. Extend $F_s$ to a path $\wt F_s$ in $\op{Aut}_{\R}(TY|_{f(\Lambda)})$. Let $R$ be a vector field in $TY|_{f(\Lambda)}$, positively transverse to $\xi$. Then $\wt F_1^{-1}(R)$ is nowhere tangent to $T\Lambda$, and the linking number of the knot with the vector field does not depend on the choice of lifting $\wt F_s$. This integer is called the \emph{Thurston-Bennequin number} of $(f, F_s)$, denoted $tb(f, F_s)$.
\end{definition}

\begin{remark} 
When $n$ is even, the definition makes sense but the invariant is uninteresting. In the example $\R^{2n+1}_{std}$, and choosing a genuine Legendrian representative, we can compute $tb$ by the signed count of self intersections in the Lagrangian projection (in any dimension). If $n$ is odd, the intersection product is skew, and the order of the inputs is given by height. For even $n$ the intersection product is commutative, so all the data necessary to calculate $tb(f)$ is contained in the Lagrangian projection. Together with the Lagrangian neighborhood theorem, it follows that $tb(f) = -\frac12\chi(\Lambda)$ in this case.
\end{remark}




\begin{prop} \label{prop:form inv}
Suppose $\Lambda$ is stably parallelizable. We describe formal Legendrian  up to formal isotopy in $\R^{2n+1}_\std$.
\begin{itemize}
\item[(a)] Suppose $n$ is odd. If two formal Legendrian embeddings have the same Thurston-Bennequin number and rotation class, then they are formally Legendrian isotopic.
\item[(b)] If two formal Legendrian surfaces in $\R^5_\std$ have the same rotation class, they are formally Legendrian isotopic. In particular all formal embeddings of $S^2$ are formally isotopic.
\item[(c)] Suppose $n>2$ is even. Then for each rotation class there are at most two formal Legendrian isotopy classes. If $\Lambda$ is simply connected, there are exactly two.
\end{itemize}
\end{prop}

\begin{remarks}
Every set of invariants is realized by a formal Legendrian embedding, with the additional note in case (a) that the parity of $tb$ is determined by $r$. However note that Corollary \ref{cor:exist} is false if $n=1$: there is no Legendrian realizing a formal Legendrian unknot with $tb = 0$. 

For $n>3$ the parity of $tb$ is determined only by the topology of $\Lambda$, for example $tb$ is odd for any Legendrian sphere in $\R^{2n+1}_{std}$. To show this, first take the Lagrangian projection of the Legendrian, which is an exact Lagrangian immersion in $\R^{2n}_\std$. Notice the parity of $tb$ is equal to the $\bmod \, 2$ count of self interesections of this Lagrangian immersion, in fact this is an invariant of \emph{smooth} immersions in $\R^{2n}$ up to regular homotopy. Both smooth and Lagrangian immersions satisfy h-principles \cite{PDR}, thus the existence of Lagrangian immersions of a given smooth regular homotopy class is governed by the inclusion map $\pi_nU_n \to \pi_nV_{2n,n}$. For $n$ odd this is a map $\Z \to \Z_2$, and (a stable shift of) Lemma \ref{class lem} implies this is the zero map except when $n=1, 3$.

It is unknown to the author if there exists a calculable invariant in $\Z_2$ which distinguishes the formal isotopy classes in case (c). Below it is defined as an invariant associated to a smooth isotopy between two Legendrian embeddings, which is why the $\pi_1\Lambda = 0$ assumption is needed. In the case of spheres, we can measure this formal class by attaching a Weinstein handle to the Legendrian, and then asking whether the resulting Weinstein manifold is diffeomorphic to $T^*S^{n+1}$ or $S^{n+1} \x \R^{n+1}$. Strangely however, when $n=6$ this formal isotopy class still exists, even though the distinction between these sphere bundles does not.
\end{remarks}

\begin{proof}[Proof of Proposition \ref{prop:form inv}:]
We assume some basic facts about frame bundles, see \cite{B, St}. We let $V_{2n+1, n}$ be the Steifel manfifold (i.e. the space of $n$--frames in $\R^{2n+1}$). We lose no generality by assuming the formal Legendrian embeddings are genuine, so we do. Given two Legendrian embeddings construct a smooth isotopy $f_t$ between them, this defines a path $\beta_t: \Lambda \to V_{2n+1,n}$ so that $\beta_0$ is a constant map. Inside $V_{2n+1, n}$, identify $U_n$ as the subset of Legendrian frames. (Though ``which frames are Legendrian'' depends on the point in $\R^{2n+1}$, these inclusions are all homotopy equvialent to the inclusion $U_n \subseteq O_{2n} \subseteq O_{2n+1} \to V_{2n+1, n}$.) $\beta_1$ has image inside of $U_n$ since $f_1$ is Legendrian, and so $\beta_t$ defines an element $\beta \in \pi_1\left(\op{Map}\left(\Lambda, V_{2n+1, n}\right), \op{Map}\left(\Lambda, U_n\right)\right)$. Notice that $\op{Map}(\Lambda, V_{2n+1, n})$ is conected since $V_{2n+1, n}$ is $n$-connected.

Our smooth isotopy can be made into a formal Legendrian isotopy exactly when $\beta = 0$. Conversely, given any $\beta \in \pi_1\left(\op{Map}\left(L, V_{2n+1, n}\right), \op{Map}\left(L, U_n\right)\right)$ and a Legendrian embedding $f_0$, we can define a formal Legendrian embedding $(f, F_s) = (f_0, \beta_s)$. If $f_1$ is a Legendrian realizing this formal Legendrian, then the obstruction associated to the given smooth isotopy between $f_0$ and $f_1$ is $\beta$.

In the long exact sequence for the pair, notice $\partial_* \beta = r(f_0) - r(f_1) \in \pi_0\op{Map}(\Lambda, U_n)$. Thus under the assumption $r(f_0) = r(f_1)$ we can lift $\beta$ to $\tilde{\beta} \in \pi_1\op{Map}(\Lambda, V_{2n+1, n})$. We pause to prove some lemmas concerning the homotopy groups of frame bundles.

\begin{lemma}\label{normal lem}
Consider the fibration $O_{n+1} \to O_{2n+1} \to V_{2n+1, n}$. In the homotopy long exact sequence, the map $\pi_{n+1}V_{2n+1, n} \to \pi_nO_{n+1}$ is injective, except for $n=2,6$. For these two values, $\pi_nO_{n+1}$ is trivial.
\end{lemma}

\begin{proof}
First, consider the case where $n$ is odd. The kernel of our map is the image of the group $\pi_{n+1}O_{2n+1}$. By Bott periodicity, this group is finite. But $\pi_{n+1}V_{2n+1, n} \cong \Z$, so the image must be trivial.

Next, consider the case where $n$ is even, and not equal to $2$ or $6$. Consider the map $\pi_nO_{n+1} \to \pi_nO_{2n+1}$. The first group classifies $(n+1)$-vector bundles on $S^{n+1}$, whereas the second group classifies stable bundles. Since $TS^{n+1}$ is non-trivial, but stably trivial \cite{BM}, we know this map must have non-zero kernel. So $\pi_{n+1}V_{2n+1, n} \to \pi_nO_{n+1}$ has non-zero image. Since $\pi_{n+1}V_{2n+1, n} \cong \Z_2$, this implies the map is injective.
\end{proof}

\begin{lemma}\label{class lem}
For all $n > 2$, $\pi_{n+1}U_n \to \pi_{n+1}V_{2n+1, n}$ is the zero map. For $n=2$, it is a surjection.
\end{lemma}

\begin{proof}
Let $n\neq2,6$. Notice that the inclusion $U_n \subseteq V_{2n+1, n}$ factors through $U_n \subseteq O_{2n+1} \to V_{2n+1, n}$. By the previous lemma, the second map is trivial on $\pi_{n+1}$.

If $n=6$, consider the map $\pi_{n+1}U_n \to \pi_{n+1}O_{2n+1}$. This is in the stable range, so we can look at the exact sequence

\[\pi_{n+1}U \to \pi_{n+1}O \to \pi_{n+1}(O/U) \to \pi_nU.\]

By Bott periodicity, $\pi_nU \cong 0$, and $\pi_{n+1}(O/U) \cong \pi_{n+1}(\Omega O) \cong \Z_2$. It follows that the map $\pi_{n+1}U_n \to \pi_{n+1}O_{2n+1}$ is multiplication by $2$, as a map $\Z \to \Z$. Therefore, the map $\pi_{n+1}U_n \to \pi_{n+1}V_{2n+1, n} \cong \Z_2$ is zero.

Case $n=2$: Since $\pi_nO_{n+1} \cong 0$, we know $\pi_{n+1}O_{2n+1}$ surjects onto $\pi_{n+1}V_{2n+1, n}$. This, together with the fact that $\pi_{n+1}U_n \to \pi_{n+1}O_{2n+1}$ is an isomorphism, implies the result.
\end{proof}

\begin{lemma}\label{tb lem}
Let $n$ be odd. From the fibrations $O_{n+1} \to O_{2n+1} \to V_{2n+1, n}$ and $O_n \to O_{n+1} \to S^n$, form the composition map $tb: \pi_{n+1}V_{2n+1, n} \to \pi_nO_{n+1} \to \pi_nS^n$. Then $tb$ is an injection, in fact, it is the map $\Z \mapsto 2\Z$.\\
\end{lemma}

\begin{proof}
We know from Lemma \ref{normal lem} that the first map is an injection, so the lemma is equivalent to the claim that $\op{Im}(\pi_{n+1}V_{2n+1, n}) \cap \ker(\to \pi_nS^n)$ is trivial in $\pi_nO_{n+1}$. By the exact sequences, this group is equal to the intersection $\ker(\to \pi_nO_{2n+1}) \cap \op{Im}(\pi_nO_n) \subseteq \pi_nO_{n+1}$. Thus the lemma is equivalent to the statement that any rank $n+1$ vector bundle on $S^{n+1}$ which is both stably trivial and zero euler class is in fact trivial. This is true, because the tangent bundle of the sphere generates the group of stably trivial vector bundles over $S^{n+1}$, and it has nonzero euler class. The second statement follows since the euler class of this generator is $2$.
\end{proof}

Returning to the proof of the theorem, recall our isotopy is unobstructed if $\tilde{\beta} \in \pi_1\op{Map}(\Lambda, V_{2n+1, n})$ is in the image of $\pi_1\op{Map}(\Lambda, U_n)$. Take any degree one map $\Lambda \to S^n$. Since $V_{2n+1, n}$ is $n$-connected this map induces an isomorphism $\pi_1\op{Map}(\Lambda, V_{2n+1, n}) \cong \pi_{n+1}V_{2n+1, n}$, identifying the image of $\pi_1\op{Map}(\Lambda, U_n)$ with that of $\pi_{n+1}U_n$.

In part (a), $n$ is odd. We claim $tb(\tilde{\beta}) = tb(f_0) - tb(f_1)$. Since Lemma \ref{tb lem} states $tb: \pi_{n+1}V_{2n+1, n} \to \pi_nS^n$ is an injection, this implies that $\tilde{\beta} = 0$ if and only if $tb(f_0) = tb(f_1)$; the desired result. Consider the geometric meaning of the maps in Lemma \ref{tb lem}. The first map to $\pi_nO_{n+1}$ can be interpreted as the difference class of the Legendrian framings of the normal bundle, with identification induced by the isotopy. The second map, induced by $O_{n+1} \to S^n$, is simply ``pick one vector in the frame'', here we think of it as choosing the Reeb vector field. Thus $tb(\tilde{\beta})$ represents the difference class of the Reeb framings, which equals $tb(f_0) - tb(f_1)$.

In part (b), $n=2$. Lemma \ref{class lem} implies that that $\tilde{\beta}$ is in the image of $\pi_{n+1}U_n$, thus $\beta = 0$.

In part (c), $n>2$ is even. $\tilde{\beta} \in \pi_{n+1}V_{2n+1, n} \cong \Z_2$, which implies there are at most two formal Legendrian isotopy classes for the given rotation class. However $\tilde{\beta}$ is an invariant of a smooth isotopy: one can imagine a isotopy from a Legendrian to itself so that $\tilde{\beta} \neq 0$. If such a case exists there will only be one formal isotopy class for the given rotation class. Under the assumption $\pi_1\Lambda = 0$, the space of smooth embeddings $\Lambda \hookrightarrow \R^{2n+1}_{std}$ is simply connected \cite{SC} and thus this cannot occur.
\end{proof}

\end{document}